\documentclass[a4paper,10pt]{article}

\usepackage{amsmath,amsfonts,amssymb,amsbsy}
\usepackage{amsthm}
\usepackage{multicol}

\numberwithin{equation}{section}

\begin{document}
\title{\bf \huge Self-organization in a group of mobile autonomous agents
\footnote{Peking University Swarm Research Group Members:
Professor Long Wang, Professor Tianguang Chu, Ms. Shumei Mu, Miss
Bo Liu and Miss Hong Shi. }}
\author{P. K. U. Swarm
\\
\small{{\textit{Center for Systems \& Control, Dept. of Mechanics
\& Engineering Science, Peking University, Beijing 100871, China}
}}}
\date{}
\maketitle

\thispagestyle{empty}

\begin{minipage}{16cm}

{\small\textbf{Abstract}}: This paper considers a discrete time
swarm model of a group of mobile autonomous agents with a simple
attraction and repulsion function for swarm aggregation and
investigates its stability properties. In particular, it is proved
that the individuals (members) of the swarm will aggregate and
form a cohesive cluster of a finite size  depending only on the
parameters of the swarm model in a finite time, and the swarm
system is completely stable.

\end{minipage}

\vspace{10mm}

\begin{multicols}{2}
\section{Introduction}

In recent years the topic of swarms has attracted considerable
attention because  swarming behavior can be found in many
organisms in nature, ranging from simple bacteria to large
animals. Generally speaking, ``swarm behavior" is a kind of
aggregate motion that a variety of organisms have the ability to
cooperatively forage for food while trying to avoid predators and
other risks. Examples include colonies of ants, flocks of birds,
and groups of animals. Understanding the operational principles of
such motions in swarms is useful in developing distributed
cooperative control, coordination, formation control, and learning
strategies for autonomous agent systems such as autonomous
multi-robots; satellite group maneuvers; multiple unmanned
undersea/aerial vehicles, etc. The increasing interest has been
motivated to study the swarming behavior and its applications, and
work on modelling of swarming behavior [1]--[4]. Specially, the
collective dynamics models have been explored in [5]--[14]. In
this paper, we study swarming behavior in terms of aggregation,
cohesion, and stability of a group of mobile autonomous agents.

In [15], Gazi and Passino proposed an ``individual-based"
continuous time model for swarm aggregation in $n$-dimensional
space with the same attraction-repulsion rule and with identical
interaction strength among these individuals (or members). They
showed that the swarm with such model has the properties of
aggregation, cohesion and stability. In this paper we develop a
discrete time swarm model of a group of mobile autonomous agents
moving in $n$-dimensional space, and analyze its aggregation and
stability properties. Particularly, we will prove that individual
agents can form a cohesive swarm with a finite size, which depends
only on parameters of the swarm model, in a finite time.

In section 2, we present the swarm model. Then we show that the
model can converge to the center and have complete stability
behavior (i.e., cohesion) in section 3 and section 4,
respectively. Finally we briefly summarize the results of this
paper in section 5.

\section{Swarm Model}

We consider a swarm of $M$ individuals (or members) in an
$n$-dimensional Euclidian space. We assume the motion is
synchronous without delays, i.e., all these individuals move
simultaneously and know the exact position of all the other
members. The model is given by
\begin{equation}
x_i(k+1)=x_i(k)+\sum_{j=1,j\neq i}^Mg(x_i(k)-x_j(k)) \label{e1}
\end{equation}
where $x_i(k)\in R^n$ represents the position of individual $i$ at
time $k$, $i\in S$,\;$k\in N$; sets $S=\{1,2,\cdots ,M\}$ and
$N=\{0,1,2,\cdots \}$; $g(\cdot )$ represents the function of
attraction and repulsion between the members. That is, the
direction and magnitude of every individual depend on the sum of
the attraction and repulsion of the other individuals on it. The
attraction-repulsion function $g(\cdot )$ is of the form
\begin{equation}
g(y(k))=-y(k)\left( a-b\exp \left( -\frac{\Vert y(k)\Vert
^2}c\right) \right) \label{e2}
\end{equation}
where $a,b,c$ are positive constants, and $\Vert y(k)\Vert
^2=y(k)^{\top }y(k)$,\;\; $\forall k\in N$.

We can see that this function is attractive for large distances
and repulsion for small distances, and the function $g(\cdot )$
changes its sign when $g(\cdot )=0$ $(b>a)$, i.e., at the set of
points as
\[
\Theta =\left\{ y=0\;\;\mathrm{or}\;\;\Vert y(k)\Vert =\sqrt{c\ln (\frac ba)}%
=\delta \right\} .
\]

\textbf{Definition 1:} The center of the swarm members is
\[
\bar{x}(k)=\frac{1}{M}\sum_{i=1}^Mx_i(k), \;\; \forall k\in N.
\]

By the symmetry property of function $g(\cdot)$ with respect to
the origin, individual $i$ moves to each other individual $j$ just
the same
amount as $j$ moves to $i$. This implies that the center $\bar{x}%
(k)$ is stationary for all $k$. Simultaneously, we have
\begin{equation}
\sum_{j=1}^M(x_i(k)-x_j(k))=M(x_i(k)-\bar{x}(k))=Me_i(k).\label{e3}
\end{equation}
{\bf Lemma 1:} The center $\bar{x}(k)$ of the swarm described in
Eqs. (\ref{e1}) and (\ref{e2}) is stationary for all $k$.

{\bf Proof:} By the symmetry of function $g(\cdot)$, we know that
$$
\sum_{i=1}^M\sum_{j=1,j\neq i}^Mg(x_i(k)-x_j(k))=0,$$ then for all
$k$. It follows that  $$\begin{array}{rl}
\bar{x}(k+1)  =& \frac{1}{M}\sum_{i=1}^Mx_i(k+1) \nonumber\\
=& \bar{x}(k) .   \end{array}$$ Hence the swarm center
$\bar{x}(k)$ is stationary.

This lemma says that the swarm described in Eqs. (\ref{e1}) and
(\ref{e2}) is not drifting on average. In the next section, we
show that the individuals will move toward the swarm center and
form a cohesive cluster around it.

\section{Swarm Aggregation}

This section presents results concerning aggregation properties of
the autonomous agents swarm modeled in Eqs. (\ref{e1}) and
(\ref{e2}).

\textbf{Definition 2:} A swarm member $i$ is a free agent if, for
all $k$,
\[
\Vert x_i(k)-x_j(k)\Vert >\delta ,\;\;\forall j\in S,\;\;j\neq i.
\]

Since the swarm center $\bar{x}(k)$ is stationary, it is noted
$\bar{x}$ for simplicity.

\textbf{Definition 3:} The error variable is defined as
\[
e_i(k)=x_i(k)-\bar{x},\;\;\forall i\in S,\;\;\forall k\in N.
\]

\textbf{Assumption 1:} All the swarm members are free agents at
any time $k$.

\textbf{Lemma 2:} If a member $i$ of the swarm is described by Eqs. (\ref{e1}%
) and (\ref{e2}), under Assumption 1, its distance to the center
$\bar{x}$ of the swarm is greater than $\delta $, i.e.,
\[
\Vert e_i(k)\Vert =\Vert x_i(k)-\bar{x}\Vert >\delta .
\]

Then, at any time $k$, it moves toward the center $\bar{x}$; in
other words, its motion is in a direction of decrease of
$\|e_i(k)\|$.

{\bf Proof:} Set
\begin{eqnarray*}
\Lambda (i,j;k) &=&\exp\left(-\frac{\Vert x_i(k)-x_j(k)\Vert ^2}c\right), \\
\Phi &=&\exp\left(-\frac{\delta ^2}{c}\right), \\
x_{ij}(k) &=&x_i(k)-x_j(k),
\end{eqnarray*}
and by (\ref{e3}) and Assumption 1, we have
\begin{eqnarray*}
x_i(k+1)=x_i(k)-\sum_{j=1,j\neq i}^Mx_{ij}(k)(a-b\Lambda
(i,j;k))\nonumber \\
\\=x_i(k)-aMe_i(k)+b\sum_{j=1,j\neq i}^Mx_{ij}(k)\Lambda (i,j;k), \nonumber\\
\\=(1-aM)e_i(k)+b\sum_{j=1,j\neq i}^Mx_{ij}(k)\Lambda
(i,j;k)+\bar{x}, \nonumber
\end{eqnarray*}
thus
\[
e_i(k+1)=(1-aM)e_i(k)+b\sum_{j=1,j\neq i}^Mx_{ij}(k)\Lambda
(i,j;k).
\]

Let $V_i(k+1)=\frac 12e_i^{\top }(k+1)e_i(k+1),$ then
\begin{eqnarray*}
{\triangle }V_i(k)&=&V_i(k+1)-V_i(k) \\
&\leq &-a\Vert e_i(k)\Vert ^2-a(M-1)\Vert e_i(k)\Vert ^2 \\
&&+\frac12(a^2+b^2\Phi^2)M^2\Vert e_i(k)\Vert ^2 \\
&&+b(1+aM)(M-1)\delta \Phi \Vert e_i(k)\Vert .
\end{eqnarray*}
So
\begin{eqnarray*}
{\triangle }V_i(k) &\leq &-a\Vert e_i(k)\Vert ^2-(M-1)\Vert
e_i(k)\Vert\\
&&\left(\left(a-\frac{(a^2+b^2\Phi ^2)M^2}{2(M-1)}\right) \Vert
e_i(k)\Vert \right. \\
&&\left.-(1+aM)b\delta \Phi \right).\end{eqnarray*} If
$$
\Vert e_i(k)\Vert  \geq \frac{1+aM}{a-\frac{(a^2+b^2\Phi ^2)M^2}{2(M-1)}}%
b\delta \Phi \geq \delta,
$$
for $M$ large enough, then
$$
{\triangle }V_i(k) \leq -a\Vert e_i(k)\Vert ^2 =-2aV_i(k).
$$
Thus the result holds.\hfill $\square $

\textbf{Remark 1:} Lemma 2 does not say that $x_i(k)$ will converge to $%
\bar{x}$ for all $i$ and $k$.

\textbf{Theorem 1:} As time progresses all the members of the
swarm described in Eqs. (\ref{e1}) and (\ref{e2})will enter into a
bounded hyperball in a finite time bound $\bar{k},$
\[
B_\varepsilon (\bar{x})=\{x_i(k):\Vert x_i(k)-\bar{x}\Vert \leq
\varepsilon \}
\]

where
\[
\varepsilon=\frac{b}{a}\sqrt{\frac{c}{2}}\exp(-\frac{1}{2}),
\]
\;\;\;\;and
\[
\bar{k}=\max_{i\in S}\left[\frac{V_i(0)}{a\varepsilon ^2}\right].
\]

{\bf Proof:} Let $V_i(k+1)=\frac 12e_i(k+1)^{\top }e_i(k+1)$, and
by Lemma 2, we know that if $ \Vert e_i(k)\Vert \geq \delta ,$

then we will have
\[
{\triangle}V_i(k)\leq 0.
\]

Note that function $\Psi =\Lambda (i,j;k)\Vert x_i(k)-x_j(k)\Vert
$ is a
bounded decreasing function of the distance with the maximum occurring at $%
\Vert x_i(k)-x_j(k)\Vert =\delta $, and we can obtain a position
independent bound by using its maximum. (i.e., solving the
continuous equation $\frac \partial {\partial y}\left( y\exp
\left( -\frac{y^2}c\right)
\right) =0$). We know that it occurs at $\Vert x_i(k)-x_j(k)\Vert =\sqrt{%
\frac c2}$. So we can obtain that ${\triangle }V_i(k)\leq 0$ as
long as
\[
\Vert e_i(k)\Vert \geq \frac ba\sqrt{\frac c2}\exp \left(- \frac
12\right) .
\]

Define $\varepsilon =\frac ba\sqrt{\frac c2}\exp \left(- \frac
12\right) $. This implies that as $k\rightarrow \infty $, $e_i(k)$
converges within the ball $B_\varepsilon \left( \bar{x}\right) $.
Notice that $i$ is arbitrary, so the result holds for all the
swarm members.

Next, the finite time bound will be presented. Since
$${\triangle }V_i(k)=V_i(k+1)-V_i(k) $$
and by ${\triangle }V_i(k)\leq -a\Vert e_i(k)\Vert ^2=-2aV_i(k)$,
we have
\begin{eqnarray*}
V_i(k+1)-V_i(0) &=&{\triangle }V_i(0)+\cdots +{\triangle }%
V_i(k) \\
&\leq &-2a(V_i(0)+V_i(1)+\cdots +V_i(k)).
\end{eqnarray*}
Since $V_i(k)=\frac 12\Vert e_i(k)\Vert ^2$, we have, for $\Vert
e_i(k)\Vert=\varepsilon ,\;\;\forall k\in N,$
\[
V_i(0)\geq 2a\cdot \frac 12\varepsilon ^2(k+1)+\frac 12\varepsilon
^2\geq ak\varepsilon ^2,
\]
thus for all $k$,
\[
k\leq \frac{V_i(0)}{a\varepsilon ^2},
\]
so we obtain the result,
\[
\bar{k}=\max_{i\in S}\left[ \frac{V_i(0)}{a\varepsilon ^2}\right].
\]\hfill $\square $

\textbf{Remark 2:} Notice that this theorem is very important not
only because it shows the aggregation of the swarm and gives an
explicit bound on the size of the swarm in a finite time bound,
but also it says that when $M$ is large enough, the bound on the
swarm size only depends on the model parameters $a,b,c$ and is
almost independent of the number of the swarm members. However, it
does not say anything about the motion of the swarm members in the
hyperball. Next we will study the issue further.

\section{Swarm Stability}

Now we prove that the swarm system described in Eqs. (\ref{e1}) and (\ref{e2}%
) will converge to its equilibrium points, and the configuration
of the swarm members converges to a constant arrangement. Now we
define the invariable set of equilibrium points of the swarm
system as
$$
\Omega _e=\left\{x:\sum_{j=1,j\neq i}^Mg(x_i-x_j)=0,\;i\in S
\right\}
$$
where $x=\{x_{1}^\top,\cdots,x_M^\top\}^\top \in R^{Mn} $.

\textbf{Theorem 2:} As $k\rightarrow \infty $, all the members of
the swarm described in Eqs. (\ref{e1}) and (\ref{e2}) converge to
$\Omega _e$.

{\bf Proof:} For all $k\in N,$ let the Lyapunov function be
$$
J(x_i(k)) =\frac 12\sum_{i=1}^{M-1}\sum_{j=i+1}^M \Vert
x_i(k)-x_j(k)\Vert ^2$$ and we have
\[
\triangle J(x_i(k))=J(x_i(k+1))-J(x_i(k)).
\]
In fact, it can be seen that $\triangle J(x_i(k))\leq 0$ as long
as
\[
\bar{\triangle}=\Vert x_i(k+1)-x_j(k+1)\Vert ^2-\Vert
x_i(k)-x_j(k)\Vert ^2\leq 0.
\]
So next we will consider $\bar{\triangle},$ and for simplicity, we
set
\begin{eqnarray*}
\theta (i,l;k) &=&\exp \left(-\frac{\Vert x_i(k)-x_l(k)\Vert ^2}c\right), \\
\phi (j,l;k) &=&\exp \left(-\frac{\Vert x_j(k)-x_l(k)\Vert
^2}c\right),
\end{eqnarray*}
\begin{eqnarray*}
x_{ij}(k)=x_i(k)-x_j(k),\\
x_{il}(k)=x_i(k)-x_l(k),\\
x_{jl}(k)=x_j(k)-x_l(k),
\end{eqnarray*}
\[
\Delta =x_i(k+1)-x_j(k+1),\;\;\Phi =\exp (-\frac{\delta ^2}c),
\]

and then
\begin{eqnarray*}
\Delta&=&(1-aM)x_{ij}(k)\\
&&+b\left( \sum_{l=1,l\neq i}^M\theta
(i,l;k)x_{il}(k)\right.\\
&&\left.-\sum_{l=1,l\neq j}^M\phi (j,l;k)x_{jl}(k)\right),
\end{eqnarray*}
so we have
\begin{eqnarray*}
\Vert \Delta \Vert ^2 &=&\Delta ^{\top }\Delta  \\
&\leq &(1-aM)^2\Vert x_{ij}(k)\Vert ^2+b^2\Phi ^2M^2\Vert
x_{ij}(k)\Vert ^2\\
&&+4b\delta (1+aM)(M-1)\Phi \Vert x_{ij}(k)\Vert.
\end{eqnarray*}
Then,
\begin{eqnarray*}
\frac 12\overline{\Delta }&=&-a\Vert x_{ij}(k)\Vert ^2-(M-1)\Vert
x_{ij}(k)\Vert  \\
&&\times \left( (a-\frac{a^2+b^2\Phi ^2}{2(M-1)}%
M^2)\Vert x_{ij}(k)\Vert \right.\\
&&\left.-2b\delta (1+aM)\Phi \right).
\end{eqnarray*}
Under Assumption 1, when $M$ is large enough, we have
$$
\Vert x_{ij}(k)\Vert \geq \delta,
$$
then
\[
\Vert \Delta \Vert ^2\leq -a\Vert x_{ij}(k)\Vert ^2\leq 0.
\]
From this,
\[
\triangle J(x_i(k))=J(x_i(k+1))-J(x_i(k))\leq 0.
\]
Thus by the LaSalle's Invariance Principle we can obtain the fact that as $%
k\rightarrow \infty $ the state $x_i(k)$ converges to the largest
invariant subset of the set defined as
\begin{eqnarray*}
\Omega  &=&\left\{x:\triangle J(x)=0\right\} \\
&=&\left\{x:\sum_{j=1,j\neq i}^Mg(x_i-x_j)=0\right\} \\
&=&\Omega _e.
\end{eqnarray*}
Because each point in $\Omega _e$ is an equilibrium, $\Omega $ is
an invariant set. This result is proved.\hfill $\square $

\textbf{Remark 3:} Notice that the stability of the swarm system
implies global convergence to its equilibrium point set. And we
know that it does not require the system to have unique
equilibrium point. Moreover, note also that in the analysis of the
above results, the dimension of the state space $n$ was not used
so that these results hold for any dimension.

\section{Conclusion}

We have shown that the swarm model described in Eqs. (\ref{e1})
and (\ref{e2}) can exhibit aggregation, cohesion, and global
convergence behavior. That is, such properties are of practical
interest in formation control of multi-robot systems. All
members of the swarm eventually enter into a finite size
ball in a finite bounded time and converge to a constant arrangement.

\end{multicols}

\end{document}